\documentclass[10pt]{amsart}
\usepackage{pgf,tikz}

\usepackage{braids}
\usetikzlibrary{decorations.pathreplacing}

\usepackage{caption}

\usepackage{amssymb}
\usepackage{stmaryrd}
\usepackage{mathrsfs}
\usepackage[all]{xy}

\makeatletter
  \def\@wrindex#1{%
    \protected@write\@indexfile{}%
      {\string\indexentry{#1}{ \S\thesubsection (p.\thepage)}}
    \endgroup
  \@esphack}

\makeatother

\makeindex

\usepackage{rotating}

\usepackage{amscd}
\usepackage{epsfig}

\author{Benjamin Enriquez}
\author{Florence Lecomte}

\address{Institut de Recherche Math\'{e}matique Avanc\'{e}e, UMR 7501, 
Universit\'{e} de Strasbourg et CNRS, 7
rue Ren\'{e} Descartes, 67000 Strasbourg, France}
\email{enriquez@math.unistra.fr}

\address{Institut de Recherche Math\'{e}matique Avanc\'{e}e, UMR 7501, 
Universit\'{e} de Strasbourg et CNRS, 7
rue Ren\'{e} Descartes, 67000 Strasbourg, France}
\email{lecomte@math.unistra.fr}

\date{\today}
\newtheorem{thm}{Théorème}[section]
\newtheorem{lem}[thm]{Lemme}

\newtheorem{prop}[thm]{Proposition}

{\theoremstyle{definition} }
{\theoremstyle{definition} \newtheorem{defn}[thm]{Définition}}

{\theoremstyle{remark} }

\numberwithin{equation}{subsection}
\numberwithin{figure}{section}
\begin{document}

\baselineskip 16pt 

\title[Sur l'interprétation homologique %du complété pro-unipotent
 du groupe fondamental]{Sur l'interprétation homologique du complété 
pro-unipotent du groupe fondamental}

\begin{abstract}
Après un rappel de l'interprétation cohomologique des quotients unipotents du groupoïde fondamental d'une variété algébrique sur un sous-corps de $\mathbb C$  (Beilinson, Deligne-Goncharov), nous proposons une construction homologique des morphismes de transition. 
%A partir de l'interprétation cohomologique  du complété pro-unipotent du groupo\-ïde fondamental d'une variété algébrique sur un sous-corps de 
%$\mathbb C$ (Beilinson, Deligne-Goncha\-rov), nous construisons une famille de morphismes de quotients de l'algèbre de ce groupoïde vers des 
%groupes d'homologie singulière relative ainsi que des morphismes entre ces groupes induits par cette structure. 
\end{abstract}

\bibliographystyle{amsalpha+}
\maketitle
%{\footnotesize \tableofcontents}

\subsection*{Introduction}

Dans le but de donner un statut motivique aux groupes fondamentaux des variétés algébriques, les auteurs de \cite{DG}, 
suivant une idée de Beilinson, ont interprété les quotients nilpotents successifs de l'algèbre du groupoïde fondamental 
d'une variété algébrique $X$ en termes du dual de l'hypercohomologie de faisceaux sur les puissances $X^n$, laquelle s'exprime 
en termes de cohomologie relative de paires $(X^n, Y^{(n)})$, où $Y^{(n)}$ est une union de sous-variétés de $X^n$. Ce travail 
a été rédigé de façon détaillée dans \cite{BGF}.    

On déduit de ces résultats une famille d'isomorphismes entre le $n$ième quotient de l'algèbre du groupoïde de $X$ et l'homologie singulière relative de $(X^n,Y^{(n)})$. Nous nous proposons de construire explicitement les opérations entre ces groupes 
d'homologie singulière relative correspondant, sous ces isomorphismes, aux opérations des quotients de l'algèbre du groupoïde 
fondamental. 

Dans cette note, nous effectuons cette construction pour les morphismes de transition entre les quotients successifs de l'algèbre du groupoïde fondamental. 

La construction est contenue dans le théorème \ref{thm:ppal}. Elle repose sur la notion de paire collier (\cite{Gr}).  Nous montrons 
que les paires $(X^n,Y^{(n)})$ sont des exemples de paires collier (proposition \ref{prop:paire:collier}). Les propriétés 
homologiques des paires collier permettent la construction du morphisme recherché entre les homologies relatives de 
$(X^{n+1},Y^{(n+1)})$ et $(X^n,Y^{(n)})$.

\subsubsection*{Remerciements} Nous remercions V. Kharlamov de nous avoir indiqué la référence \cite{L} et expliqué son 
contenu. B.E. est membre du projet ANR HighAGT ANR-20-CE40-0016. 

%Sujet : expression des quotients nilpotents du groupoïde fondamental en termes homologiques. 
%
%Existant : hypercohomologies de faisceaux identifiées à cohomologie relative. 
%
%Projet : interprétation de BLABLA en termes d'homologie singulère relative. Premier pas : applications 
%du groupoide fondamental vers groupes d'homologie singulère relative qui sont isos grâce à littérature. 
%
%Traductions en termes de groupes d'homologie singulère relative de diverses opérations sur groupoïdes.  
%
% Nos projets pour la suite. 

\subsection{L'expression des quotients nilpotents du groupoïde fondamental en termes homologiques}
\label{sect:0:1}

Soit $X$ une variété algébrique lisse complexe connexe et $x,y\in X$. Pour $n\geq1$, 
on note $Y^{(n)}\subset X^n$ la sous-variété donnée par $Y^{(n)}:=Y_0^{(n)}\cup\cdots\cup Y_n^{(n)}$, avec $Y_i^{(n)}\subset X^n$ définie par 
l'équation $x_i=x_{i+1}$ (avec $x_0:=x$, $x_{n+1}=:y$). 

Le résultat suivant, d\^u à Beilinson, a été rédigé dans \cite{DG}, puis en plus de détail dans 
\cite{BGF}. 

\begin{thm}\label{thm:bgf} (\cite{BGF}, Thm. 3.184 et p. 207 et \cite{DG}, Prop. 3.4)
On a un morphisme de $\mathbb Q$-espaces vectoriels
\begin{equation}\label{appl:ppale}
\sigma_{yx} :H^n(X^n,Y^{(n)};\mathbb Q)\to(\mathbb Q\pi_1(X;y,x)/\mathbb Q\pi_1(X;y,x)\cdot J^{n+1})^\vee
\end{equation}
où $\pi_1(X;y,x)$ est l'ensemble des classes d'homotopie de chemin tracés sur $X$ avec début en $x$ et fin en $y$, 
où $J:=\mathbb Q\pi_1(X;x)_+$ est l'idéal d'augmentation de l'algèbre de groupe de $\pi_1(X;x):=\pi_1(X;x,x)$, où 
$\cdot$ est l'action de $\mathbf k\pi_1(X;x)$ sur $\mathbf k\pi_1(X;y,x)$ par composition des chemins, et où 
$H^n(-,-;\mathbb Q)$ est la cohomologie singulière relative. 

Si $x\neq y$, alors $\sigma_{yx}$ est un isomorphisme. Si $x=y$, on a une suite exacte 
$$
0\to H^n(X^n,Y^{(n)};\mathbb Q)\to(\mathbb Q\pi_1(X;x)/J^{n+1})^\vee\to\mathbb Q\to0.
$$ 
\end{thm}

\begin{defn}
On note 
$$
\tau_{yx}^{(n)}:\mathbb Q\pi_1(X;y,x)/\mathbb Q\pi_1(X;y,x)J^{n+1}\to H_n(X^n,Y^{(n)};\mathbb Q)
$$
le morphisme de $\mathbb Q$-espaces vectoriels donné par la composition de l'application canonique de 
$\mathbb Q\pi_1(X;y,x)/\mathbb Q\pi_1(X;y,x)J^{n+1}$ dans son double dual, du morphisme de $\mathbb Q$-espaces 
vectoriels $\sigma_{yx}^\vee$ dual de $\sigma_{yx}$, et de l'isomorphisme $H^n(X^n,Y^{(n)};\mathbb Q)^\vee\simeq H_n(X^n,Y^{(n)};\mathbb Q)$ (voir \cite{Gr}, (23.28)). 
\end{defn}

\begin{lem}
Si $x\neq y$, alors $\tau_{yx}^{(n)}$ est un isomorphisme de $\mathbb Q$-espaces vectoriels. Si $x=y$, on a une suite exacte 
$0\to\mathbb Q\to \mathbb Q\pi_1(X;x)/J^{n+1}\stackrel{\tau_{xx}^{(n)}}{\to} H_n(X^n,Y^{(n)};\mathbb Q)\to 0$, où la 
première application envoie $1$ sur la classe de l'unité de $\pi_1(X;x)$.  
\end{lem}

\proof Il résulte du théorème \ref{thm:bgf} que la source de $\tau_{yx}^{(n)}$ est de dimension finie, ce qui implique qu'elle est isomorphe à son double dual, et donc que $\tau_{yx}^{(n)}$ a les propriétés annoncées. \hfill\qed\medskip

Si $\gamma$ est un chemin tracé sur $X$ avec début en $x$ et fin en $y$, alors $\gamma^{(n)}:\Delta_n\to X^n$ donné par 
$\Delta_n=\{(t_1,\ldots,t_n)|0\leq t_1\leq\ldots\leq t_n\leq 1\}\ni (t_1,\ldots,t_n)\mapsto (\gamma(t_1),\ldots,\gamma(t_n))\in X^n$ 
est un cycle relatif de $(X^n,Y^{(n)})$, les images des faces $t_i=t_{i+1}$, $0=t_0$, $t_n=0$ étant contenues dans $Y^{(n)}$. 
De plus, le morphisme $\tau_{yx}^{(n)}$ envoie la classe de $\gamma\in\pi_1(X;y,x)$ vers la classe de $\gamma^{(n)}$ dans 
l'homologie relative $H_n(X^n,Y^{(n)};\mathbb Q)$. 

\subsection{Paires collier : propriétés et exemples}

%\subsubsection{Définitions}

Soit $(M,A)$ un couple formé d'un espace topologique $M$ et d'un sous-espace $A$. 

\begin{defn} \label{def:paire:collier}
([Gr], \S19, p. 112). $(M,A)$ est une {\it paire collier} si et seulement si on a

(a) $M$ est séparé ; 

(b) les points de $M-A$ peuvent être séparés de $A$ : pour tout $x\in M-A$, il existe des voisinages $U$ de $x$ et $V$ de $A$ 
qui sont disjoints ; en particulier, $A$ est fermé dans $M$ ; 

(c) il existe un voisinage ouvert $C$ de $A$ dans $M$ contenant strictement $A$, 
tel que  $A$ est une rétraction forte par déformation de $C$, \`a savoir il existe une 
application $r: C\rightarrow A$ telle que $i\circ r$ est l'identité de $A$, où $i$ est l'inclusion $i: A\rightarrow C$, 
et $i\circ r$ est homotope \`a l'identit\'e de $C$ relativement \`a  $A$.
\end{defn}

\begin{prop}\label{prop:paire:collier}
Si $M$ est une vari\'et\'e analytique complexe et $A$ est une union finie de sous-variét\'es analytiques, alors $(M,A)$ est une 
paire collier.% (cf. définition \ref{def:paire:collier}).   
%Pour tout entier $n\geq 1$, la paire 
%\begin{equation}\label{PC:n}
%(X^n , Y^{(n)})
%\end{equation} 
%est une paire collier . 
\end{prop}

\proof Les points (1)-(2) sont clairs. Montrons le point (3). D'apr\`es \cite{L}, Thm. 2,  on peut trianguler la paire
$(M,A)$, qui est donc une paire polyédrale au sens de \cite{S}, p. 113. D'après \cite{S}, \S3.3, Cor. 11, ceci implique que 
$A$ est une rétraction forte par déformation d'un de ses voisinages dans $M$. \hfill\qed\medskip 

%[Réservoir d'exemples $(X^n$,une union des $X^{n-1})$ p ex $(X^n , Y^{(n)})$ avec $n\geq1$]

%\subsubsection{Propriétés des paires collier}

\begin{thm} \label{thm:excision} (cf. \cite{Gr}, Prop. 19.3) 
(a) Soit $(M,A)$ une paire collier et $f:A\to B$ un morphisme d'espaces topologiques avec $Y$ séparé. 
Soit $N:=M\cup_f B$. Alors $(N,B)$ est une paire collier. 

(b) (cf. \cite{Gr}, Thm. 19.14) De plus $(M,A)\to(N,B)$ induit un isomorphisme $H_*(M,A;\mathbb Q)\stackrel{\sim}{\to}
H_*(N,B;\mathbb Q)$ 
entre les homologies relatives. 
\end{thm}

%[Mayer-Vietoris (ds Spanier-ds Greenberg)] 

\subsection{Les versions homologiques des morphismes de transition entre quotients du groupoïde fondamental}

Soit $n\geq1$.  Avec les notations de \S\ref{sect:0:1}, il suit de la proposition \ref{prop:paire:collier} que 
\begin{equation}\label{PC:n}
(X^{n-1}, Y^{(n-1)})
\end{equation} 
est une paire collier.  

%La construction de la fl\`eche entre les homologies relatives repose essentiellement sur : 
\begin{thm} \label{thm:ppal}
Soit $n\geq1$.  La paire collier \eqref{PC:n} donne naissance à un morphisme
$$
\kappa_n: H_n(X^n, Y^{(n)};\mathbb Q) \longrightarrow H_{n-1}(X^{n-1} , Y^{(n-1)};\mathbb Q)
$$
qui s'insère dans le diagramme
commutatif suivant
\begin{equation}\label{CD}
\xymatrix{ \mathbb Q\pi_1(X;y,x)/\mathbb Q\pi_1(X;y,x)J^{n+1}\ar^{\ \ \ \ \ \ \ \ \tau_{yx}^{(n)}}[r]\ar[d]&  
H_{n}(X^{n},Y^{(n)};\mathbb Q)\ar^{\kappa_n}[d]\\ 
\mathbb Q\pi_1(X;y,x)/\mathbb Q\pi_1(X;y,x)J^{n} \ar_{\ \ \ \ \ \ \ \ \ \ \tau_{yx}^{(n-1)}}[r]&  H_{n-1}(X^{n-1},Y^{(n-1)};
\mathbb Q)}
\end{equation} 
\end{thm}

\proof 
On pose $Z^{(n)} = Y_0^{(n)}\cup Y_1^{(n)}\cup \dots \cup Y_{n-1}^{(n)}\hookrightarrow Y^{(n)}$ et 
$A= Y_{n}^{(n)}\cap Z^{(n)}$. 

L'isomorphisme $X^{n-1}\simeq X^{n-1}\times\{y\}=Y_{n}^{(n)}$ induit un isomorphisme 
$(X^{n-1},Y^{(n-1)})\simeq (Y_{n}^{(n)},A)$, ce qui implique par la proposition 
\ref{prop:paire:collier} que $(Y_{n}^{(n)},A)$ est une paire collier.
D'apr\`es le théorème \ref{thm:excision} (a), il suit que l'extension 
$(Y_{n}^{(n)}\cup Z^{(n)},Z^{(n)}) = (Y^{(n)}, Z^{(n)})$ est \'egalement 
une paire collier et de plus par le théorème \ref{thm:excision} (b), le morphisme induit
$$ 
H_q (Y_{n}^{(n)},A;\mathbb Q) \rightarrow H_q (Y^{(n)}, Z^{(n)};\mathbb Q)
$$
 est un isomorphisme pour tout entier $q$, ce qui en spécialisant $q$ et en utilisant 
l'isomorphisme $(X^{n-1},Y^{(n-1)})\simeq (Y_{n}^{(n)},A)$ fournit un isomorphisme 
 $$
g_{n-1}:H_{n-1}(X^{n-1} , Y^{(n-1)};\mathbb Q) \simeq H_{n-1} (Y^{(n)}, Z^{(n)};\mathbb Q)  .
$$

En composant l'inverse de cet isomorphisme avec le morphisme connectant 
$H_{n}(X^{n},Y^{(n)};\mathbb Q)\to H_{n-1}(Y^{(n)},Z^{(n)};\mathbb Q)$ 
associé au triplet $Z^{(n)}\subset Y^{(n)}\subset X^{n}$,
on obtient le morphisme recherché 
$$
H_{n}(X^{n},Y^{(n)};\mathbb Q)\to H_{n-1}(X^{n-1} , Y^{(n-1)};\mathbb Q).
$$
 
%On part du morphisme connectant 
%$$
%f_n: H_{n+1}(X^{n+1} , Y^{(n+1)}) \rightarrow H_n(Y^{(n+1)})
%$$ 
%et on pose $Z = Y_1^{(n+1)}\cup Y_2^{(n+1)}\cup \dots \cup Y_{n+1}^{(n+1)}\hookrightarrow Y^{(n+1)}$ et $A= Y_0^{(n+1)}\cap Z$. 
%
%Comme on a $Y_0^{(n+1)}\simeq X^n$ et $A\simeq Y^{(n)}$ , $(Y_0^{(n+1)}, A)$ est, par la proposition \ref{prop:paire:collier}, 
%une paire collier. D'apr\`es \cite{Gr}, (19.3), l'extension 
%$(Y_0^{(n+1)}\cup Z,Z) = (Y^{(n+1)}, Z)$ est \'egalement une paire collier et de plus (loc. cit. 19.14 p 116) le morphisme induit
% $$ H_q (Y_0^{(n+1)},A) \rightarrow H_q (Y^{(n+1)}, Z)$$
% est un isomorphisme pour tout entier $q$, ce qui fournit un isomorphisme 
% $$g_n:H_n(X^n , Y^{(n)}) \simeq H_n (Y^{(n+1)}, Z)  .$$
% 
%Soit $h_n: H_n(Y^{(n+1)}) \rightarrow  H_n (Y^{(n+1)}, Z)$ la fl\`eche canonique. Le morphisme recherch\'e $\kappa_n$ 
%est la compos\'ee $g_n^{-1}\circ h_n \circ f_n$.

Soit $\gamma$ un chemin  de $\mathbb Q\pi_1(X;y,x)$ et $[\gamma]_{n+1}$ sa classe dans $\mathbb Q\pi_1(X;y,x)/\mathbb Q\pi_1(X;y,x)J^{n+1}$. L'image de  $[\gamma]_{n+1}$ par $\tau_{yx}^{(n)}$ dans 
 $H_{n}(X^{n} , Y^{(n)};\mathbb Q)$ est la classe de la chaîne singulière 
 \begin{align*}
 & \gamma^{(n)}:\Delta_{n}\rightarrow  X^{n},  
 \\ & (t_1,\ldots,t_{n})\mapsto (\gamma(t_1),\ldots,\gamma(t_{n}))
 \end{align*}
 qui s'envoie par $f_{n-1}$ sur son bord
 $$
 \sum_{i=0}^{i=n} (-1)^i \gamma^{(n)} \circ \partial_i, 
$$
où pour $i\in\{0,\ldots,n-1\}$, $\partial_i:\Delta_{n-1}\to\Delta_{n}$ est l'application qui envoie $(t_1,\ldots,t_{n-1})$ sur 
$(t_1,\ldots,t_i,t_i,\ldots,t_{n-1})$
pour $i\neq 0,n-1$, sur $(0,t_1,\ldots,t_{n-1})$ pour $i=0$, et sur $(t_1,\ldots,t_{n-1},1)$ pour $i=n-1$. 

 Par construction, l'application 
 $$
 \gamma^{(n)} \circ \partial_i: \Delta_{n-1}  \rightarrow  X^{n}, 
 $$
donnée par 
$$
(t_1,\ldots,t_{n-1})\mapsto \left\{\begin{array}{ll}
 (\gamma(t_1),\ldots,\gamma(t_i),\gamma(t_i)\ldots\gamma(t_{n-1}))  & \mbox{si }i\neq 0,n-1, \\
   (x,\gamma(t_1),\ldots\gamma(t_{n-1}))& \mbox{si }i=0, \\  (\gamma(t_1),\ldots\gamma(t_{n-1}),y)& \text{ si }i=n-1,
    \end{array}
\right. 
$$
s'identifie à la composition de $appl_i: X^{n-1}\simeq Y_i^{(n)} \hookrightarrow X^{n}$ et de $\gamma^{(n-1)}: \Delta_{n-1}  \rightarrow  X^{n-1}$.

Ainsi, l'image de $[\gamma]_{n+1}$ dans $H_{n-1}(Y^{(n)};\mathbb Q)= H_{n-1}(X^{n-1} \cup Z;\mathbb Q)$ est la classe de 
 $$\sum_{i=0}^{i=n} (-1)^i appl_i\circ \gamma^{(n-1)}$$
qui s'envoie par $H_{n-1}(X^{n-1} \cup Z;\mathbb Q)\to H_{n-1} (X^{n-1} \cup Z, Z;\mathbb Q)$ sur la classe de 
 $$
 appl_0 \circ \gamma^{(n-1)}: \Delta_{n-1} \rightarrow X^{n-1} \hookrightarrow X^{n-1} \cup Z 
 $$
 compte tenu de $Z= Y_1^{(n)}\cup\cdots\cup Y_{n}^{(n)}$. 
  
L'image de cette dernière classe par $g_{n-1}^{-1}$ est alors la classe de $\gamma^{(n-1)}$ dans 
$H_{n-1} (X^{n-1}, X^{n-1} \cap Z;\mathbb Q) 
= H_{n-1} (X^{n-1}, Y^{(n-1)};\mathbb Q)$, c'est-\`a-dire celle de  $[\gamma]_{n}$ par  $\tau_{yx}^{(n-1)}$.
Ceci montre la commutativité du diagramme \eqref{CD}. 
\hfill\qed\medskip 

%\subsection{Notes de travail}
%
%On note $Y_{ab}:=aX+X_\Delta+Xb\subset X^2$. On a une application 
%$$
%\kappa'_2 : H_2(X^2,Y_{ab})\to H_1(X,a+b),  
%$$
%induite par suite d'applications
%$$
%\kappa'_2 : H_2(X^2,Y_{ab})\to H_1(Y_{ab},aX+X_\Delta)=H_1(Xb,ab+bb)=H_1(X,a+b),  
%$$
%où la première application provient du triple $(X^2,Y_{ab},aX+X_\Delta)$ et la deuxième application provient  
%de l'excision dans la paire collier $(Y_{ab},aX+X_\Delta)$. 
%%, enfin la dernière égalité (non-canonique) provient de la suite exacte
%%$0=H_1(a+b)\to H_1(X)\to H_1(X,a+b)\to H_0(a+b)\to H_0(X)\to H_0(X,a+b)$ et de l'identificiation de   
%%$H_0(a+b)\to H_0(X)$ avec l'application somme $\mathbb Z\oplus\mathbb Z\to\mathbb Z$, ce qui induit à la fois 
%%$H_0(X,a+b)=0$ et une suite exacte $0\to H_1(X)\to X_1(X,a+b)\to\mathbb Z\to 0$. 
%
%On définit également l'application 
%$$
%\kappa''_2 : H_2(X^2,Y_{ab})\to H_1(X,a+b),  
%$$
%induite par la suite
%$$
%\kappa''_2 : H_2(X^2,Y_{ab})\to H_1(Y_{ab},X_\Delta+Xb)=H_1(aX,aa+ab)=H_1(X,a+b),  
%$$
%où la première application provient du triple $(X^2,Y_{ab},X_\Delta+Xb)$, et la deuxième application provient  
%de l'excision dans la paire collier $(Y_{ab},X_\Delta+Xb)$. 
%%, enfin la dernière égalité (non-canonique) est celle construite plus haut. 
%
%Pour comparer ces applications, on considère la suite exacte associée au triple $(X^2,Y_{ab},X_\Delta)$ : 
%$$
%H_2(X^2,Y_{ab})\to H_1(Y_{ab},X_\Delta)\to H_1(X^2,X_\Delta). 
%$$ 
%On a par paire collier $H_1(Y_{ab},X_\Delta)=H_1(aX+Xb,aa+bb)$, par ailleurs 
%$H_1(X^2,X_\Delta)=H_1(X)$. 
%
%
%
%
%
%
%
%
%%***
%%
%%Hopf si possible. 

\end{document}